# A modified BFGS quasi-Newton iterative formula

W. Chen

**Present mail address** (as a JSPS Postdoctoral Research Fellow): Apt.4, West 1$^{st}$ floor, Himawari-so, 316-2, Wakasato-kitaichi, Nagano-city, Nagano-ken, 380-0926, JAPAN

Permanent affiliation and mail address: Dr. Wen CHEN, P. O. Box 2-19-201, Jiangshu University of Science & Technology, Zhenjiang City, Jiangsu Province 212013, P. R. China

Present e-mail: chenw@homer.shinshu-u.ac.jp

Permanent email: chenwwhy@hotmail.com

**Abstract**

The quasi-Newton equation is the very basis of a variety of the quasi-Newton methods. By using a relationship formula between nonlinear polynomial equations and the corresponding Jacobian matrix. presented recently by the present author, we established an exact alternative of the approximate quasi-Newton equation and consequently derived an modified BFGS updating formulas.

## 1. Introduction

Recently, the present author [1,2] proved a relationship theorem as stated below

**Theorem 1.** If $N^{(m)}(x)$ and $J^{(m)}(x)$ are defined as m-order nonlinear polynomial function vector and its corresponding Jacobian matrix, respectively, then $N^{(m)}(x) = \frac{1}{m} J^{(m)}(x) x$ is always satisfied.

In this paper, we will apply this theorem to establish an exact alternative of the traditional approximate quasi-Newton equation. It is well known that the quasi-Newton equation is the very basis of various quasi-Newton methods [3,4]. Therefore, by using the presented alternative equation, we derive the modified BFGS quasi-Newton updating formulas. It is noted that the present BFGS formulas are also different from the ones previously given in Chen [1].

## 2. Modified BFGS quasi-Newton formula

To avoid time-consuming evaluation and inversion of the Jacobian matrix in each iterative step of the standard Newton method, the quasi-Newton method was developed with maintaining a superlinear convergence rate. This key of such methods is a matrix-updating procedure, of which the BFGS method is the most successful and widely used. The so-called quasi-Newton equation is the very fundamental of various quasi-Newton methods, namely,

$$J_k(x_k - x_{k-1}) = f(x_k) - f(x_{k-1}). \tag{1}$$

The Jacobian matrix $J_k$ is updated by adding a rank-one matrix to the previous $J_{k-1}$ in satisfying equation (1) and the following relations:

$$J_i r = J_{i-1} r, \quad \text{when } (x_i - x_{i-1})^T r = 0, \tag{2}$$

where $x_k - x_{k-1} = p$, $f(x_k) - f(x_{k-1}) = q$. It is emphasized that J here is the Jacobian matrix of total system. It is noted that equation (1) is an approximate one. For the polynomial-only problems, we can gain an exact alternative of equation (1) by using theorem 1. Without loss of generality, Let us consider the following nonlinear polynomial equations

$$f(x) = Lx + N^{(2)}(x) + N^{(3)}(x) + b = 0, \tag{3}$$

where $Lx$, $N^{(2)}(x)$ and $N^{(3)}(x)$ represent the linear, quadratic and cubic terms of the system of equations. The Jacobian matrix of the system is given by

$$J = \frac{\partial f(x)}{\partial x} = L + \frac{\partial N^{(2)}(x)}{\partial x} + \frac{\partial N^{(3)}(x)}{\partial x}. \tag{4}$$

By using theorem 1, we have

$$Jx = Lx + 2N^{(2)}(x) + 3N^{(3)}(x) = \bar{f}(x). \tag{5}$$

Therefore, we can exactly establish

$$J_k x_k - J_{k-1} x_{k-1} = \bar{f}(x_k) - \bar{f}(x_{k-1}) = y \tag{6}$$

After some simple deductions, we get

$$J_k(x_k - x_{k-1}) = -J_k x_{k-1} + J_{k-1} x_{k-1} + y. \tag{7}$$

It is worth stressing that equation (7) differs from equation (1) in that it is exactly constructed and forms the basis of this paper. By analogy with the deduction of the BFGS formula, let

$$J_k = J_{k-1} + u_k p_k^T. \tag{8}$$

Substituting the above equation (8) into equation (7) yields

$$u_k = \frac{1}{p_k^T x_k}(y - J_{k-1}p_k). \tag{9}$$

Therefore, we have

$$J_k = J_{k-1} + \frac{1}{p_k^T x_k}(y - J_{k-1}p_k)p_k^T. \tag{10}$$

Note that the left second term of the above equation (10) is a rank-one matrix. By using the known Shermann-Morrison formula, we can derive

$$J_k^{-1} = J_{k-1}^{-1} - \frac{(J_{k-1}^{-1}s_k)(p_k^T J_{k-1}^{-1})}{1 + p_k^T J_{k-1}^{-1} s_k}, \tag{11}$$

where $s_k = \frac{1}{p_k^T x_k}(y - J_{k-1}p_k)$. Note that y can be evaluated easily without $J_k$. The updating formulas (10) and (11) are a modified version of the following original BFGS formulas:

$$J_k = J_{k-1} - \frac{(J_{k-1}p_k - q_k)p_k^T}{p_k^T p_k} \tag{12}$$

and

$$J_k^{-1} = J_{k-1}^{-1} - \frac{(J_{k-1}^{-1}q_k - p_k)p_k^T J_{k-1}^{-1}}{p_k^T J_{k-1}^{-1} q_k}, \tag{13}$$

## 3. Remarks

One can find that the computing effort of both Eqs. (11) and (13) are nearly the same, only about $3n^2$ operations. Theoretically, the present updating formulas improve the accuracy of the solution by establishing themselves on the exact equation (7) instead of the approximate quasi-Newton equation (1). The basic idea of the quasi-Newton method is a successive update of rank one or two. Therefore, it is noted that equation (7) is not actually exact due to the approximate Jacobian matrix yielded in the previous iterative steps. It may be better to initialize Jacobian matrix J via an exact approach. In addition, it is a puzzle for a long time why one-rank BFGS updating formula performs much better compared to the other two-rank updating schemes such as DFP method [5]. In our understanding, the most possible culprit is due to the inexactness of quasi-Newton equation (1). Therefore, this suggests that the updating formulas of higher order based on exact equation (7) may be more attractive, which will include more additional curvature information to accelerate convergence. It is noted that in one dimension, the present equations (10) and (11) degenerates into the original Newton method by comparing with the fact that the traditional quasi-Newton method becomes the secant method. The performances of the present methodology need be examined in solving various benchmark problems.

Finally, it is stressed again that the present modified BFGS updating formulas (10) and (11) are different from the previously given ones in Chen [1], although they were derived based on the same equation (7). Also, the present BFGS formulas are only applicable to the nonlinear polynomial equations.